\documentclass[a4paper, oneside, 11pt]{article}
\usepackage[english]{babel}
\usepackage{amsfonts}
\usepackage{amssymb}
\usepackage{graphics}
\usepackage{amsrefs}
\usepackage{latexsym}
\begin{document}

\title{{\bf{\Large{Conditional $\alpha$-diversity for exchangeable Gibbs partitions driven by the stable subordinator}}}\footnote{{\it AMS (2000) subject classification}. Primary: 60G58. Secondary: 60G09.} }
\author{\textsc {Annalisa Cerquetti}\footnote{Corresponding author, SAPIENZA University of Rome, Via del Castro Laurenziano, 9, 00161 Rome, Italy. E-mail: {\tt annalisa.cerquetti@gmail.com}}\\
\it{\small Department of Methods and Models for Economics, Territory and Finance}\\
  \it{\small Sapienza University of Rome, Italy }}
\newtheorem{teo}{Theorem}
\date{\today}
\maketitle{}

\begin{abstract}
Asymptotic behaviour of conditional $\alpha$ diversity for the two-parameter Poisson-Dirichlet partition model and for the normalized generalized Gamma model has been recently investigated in Favaro et al. (2009, 2011) with a view to possible applications in Bayesian treatment of species richness estimation. Here we generalize those results to the larger class of mixed Poisson-Kingman species sampling models driven by the stable subordinator (Pitman, 2003).
\end{abstract}

\section{Introduction}

The concept of {\it $\alpha$-diversity} for infinite exchangeable random partitions $\Pi$ of the positive integers, induced by randomly sampling from an almost surely discrete probability measure $P$, was first introduced in Pitman (2003, cfr. Sect. 6.1 Prop. 13) as the random variable $S$, with $0 < S < \infty$,  such that, almost surely for $n \rightarrow \infty$,
$$
\frac{K_n}{n^{\alpha}} {\longrightarrow} S
$$
for $K_n$ the number of blocks in a random partition of $\Pi_n$ of $[n]$ induced by $\Pi$ and $\alpha \in (0, 1)$. For infinite random partitions induced by a r.p.m $P$ whose ranked atoms follow a Poisson-Kingman distribution $PK(\rho_\alpha, \gamma)$ driven by the L\'evy density  of the stable subordinator $\rho_\alpha= \alpha \Gamma(1-\alpha)^{-1}x^{-\alpha -1}dx$, for $0 <\alpha <1$, with $\gamma$ on $(0,\infty)$ some mixing  density, Pitman shows that $S=T^{-\alpha}$ where $T=S^{-1/\alpha}$ is the random total sum of the ranked atoms of $P$ and has density $\gamma$. (See Pitman, 2006, for a comprehensive reference on exchangeable random partitions).  Recently interest in {\it conditional $\alpha$-diversity} has emerged in  posterior species richness estimation in a Bayesian nonparametric approach to species sampling problems (cfr. Lijoi et al., 2007, 2008; Cerquetti, 2009). Given $K_n=k$ the number of species in a partition  $(n_1, \dots, n_k)$ induced by a {\it basic} sample of observed species, {\it conditional $\alpha$-diversity} is defined as the random variable $S_{\alpha}^{n,k}$ such that, almost surely for $m \rightarrow \infty$ 
$$
\frac{K_m}{m^{\alpha}} \Big\vert (K_n=k)  \rightarrow S_{\alpha}^{n, k}
$$
for $K_m$ the unknown number of {\it new} species induced by an additional sample of observations of dimension $m$.

In particular Favaro et al. (2009, 2011) derive distributional results respectively for the conditional $\alpha$ diversity under {\it two-parameter $(\alpha, \theta)$ Poisson-Dirichlet priors} (Pitman and Yor, 1997), which are well known to correspond to r.p.m.s whose ranked atoms follow a $PK(\rho_\alpha, \gamma_{\alpha, \theta})$ distribution for $\gamma_{\alpha, \theta}(t)= \frac{\Gamma(\theta +1)}{\Gamma(\theta/\alpha +1)} t^{-\theta}f_\alpha(t)$ for $0 <\alpha <1$ and $\theta > -\alpha$, and under {\it normalized generalized Gamma priors}, whose ranked atoms follow a  Poisson-Kingman law with mixing density the exponentially tilted version of the $\alpha$-stable density (cfr. Cerquetti, 2007; Lijoi et al. 2007a). 
%   and a scale mixture representation of the law of $S_{\alpha, \theta}^{n,k}$ which allows to obtain {\it highest posterior density} intervals, a Bayesian measure of uncertainty of point estimates. 
In Cerquetti (2011) a first alternative derivation for the conditional $\alpha$ diversity of the two-parameter Poisson-Dirichlet family has been obtained by a decomposition approach properly exploiting a characterization of those models in terms of the {\it deletion of classes} property (cfr. Pitman, 2003; Gnedin et al. (2009)). Here we obtain the general result for the entire class of mixed Poisson-Kingman models driven by the stable subordinator, for an arbitrary mixing density $\gamma$ on $(0, \infty)$ written as, without loss of generality, $\gamma(t)=h(t)f_\alpha(t)$ for some non negative function $h$ such that $\gamma$ is a proper density and $f_\alpha$ the density of the stable subordinator.  This class includes infinitely many priors, and as from Gnedin and Pitman (2006), induces exchangeable random partitions in Gibbs form of type $\alpha$ i.e. uniquely characterized by exchangeable partition probability function (EPPF) in the following Gibbs product form  
\begin{equation}
\label{eppf}
p(n_1,\dots, n_k)= V_{n,k} \prod_{i=1}^k (1-\alpha)_{n_j},
\end{equation}
where $(a)_b=a(a+1)\cdots(a+b-1)$ are rising factorials and the $V_{n,k}$ are weights satisying the backward recursive relation $V_{n,k}= (n -k\alpha)V_{n+1, k}+V_{n+1, k+1}$. 

\section {Main results}

The basic result in Gnedin and Pitman (2006, cfr Th. 12) states that the EPPF of each exchangeable Gibbs partition of type $\alpha$, for $\alpha \in (0,1)$, arises by mixing the EPPF corresponding to the Poisson-Kingman $PK(\rho_\alpha |T=t)$ model with the density $\gamma$ identifying the specific $PK(\rho_\alpha, \gamma)$ family, hence
$$
p_{\alpha, \gamma}(n_1, \dots, n_k)= \int_0^\infty p_{\alpha}(n_1, \dots, n_k| t) \gamma(t)dt.
$$ 
The explicit formula for $p_{\alpha}(n_1, \dots, n_k|t)$ was first obtained in Pitman (2003, cfr. eq. (66)) and is central to our main result.\\\\
{\bf Theorem 1.} Let $\Pi$ be a $PK(\rho_\alpha, \gamma)$ partition of $\mathbb{N}$ for some  $0 < \alpha <1$ and some mixing probability distribution $\gamma$ on $(0, \infty)$. Without loss of generality assume $\gamma(t)=h(t)f_\alpha(t)$. Fix $n \geq 1$ and a partition $(n_1, \dots, n_k)$ of $[n]$ with $k$ positive box-sizes, then  $\Pi_{}$ has conditional $\alpha$-diversity $S_{\alpha,h}^{n,k}$ with density 
%=T_{n,k}^{-\alpha}$ where $T_{n, k}$ has distribution
\begin{equation}
\label{conddiv}
f_{n,k}^{h,\alpha}(s)= \frac{h(s^{-1/\alpha}) \tilde{g}_{n,k}^\alpha(s)}{\mathbb{E}_{n,k}^{\alpha}[h(S^{-1/\alpha})]}
\end{equation}
for  
\begin{equation}
\label{gnka}
\tilde{g}_{n,k}^{\alpha}(s)= \frac{\Gamma(n)}{\Gamma(n-k\alpha) \Gamma(k)} s^{k-1/\alpha -1} \int_0^1 p^{n-1-k\alpha} f_\alpha((1-p)s^{-1/\alpha})dp 
\end{equation} the density of the product of independent r.v.s $Y_{\alpha, k} \times [W]^\alpha$, where  $Y_{\alpha, k}$  has density 
\begin{equation}
\label{polmit}
g_{\alpha, k\alpha}(y)=\frac{\Gamma(k\alpha +1)}{\Gamma(k+1)}y^kg_{\alpha}(y)
\end{equation}
for $g_\alpha(y)=\alpha^{-1}y^{-1-1/\alpha}f_\alpha(y^{-1/\alpha})$, and $W \sim \beta(k\alpha, n-k\alpha)$.\\\\
{\it Proof}:
By the unconditional result in Pitman (2003) recalled in the Introduction, the $\alpha$-diversity for a general $\gamma(t)=h(t)f_\alpha(t)$ mixing density, is the r.v. $S_{\alpha, h}$ with density 
\begin{equation}
\label{prior}
\gamma(s^{-1/\alpha})=h(s^{-1/\alpha})f_\alpha(s^{-1/\alpha})\alpha^{-1}s^{-1/\alpha -1}.
\end{equation}
As from Pitman (2003, cfr. Prop. 9), given $S_\alpha=s$, the exchangeable partition probability function for a $PK(\rho_\alpha|s)$ model corresponds to 
$$
p_\alpha(n_1, \dots, n_k| s^{-1/\alpha})= \frac{\alpha^k}{\Gamma(n -k\alpha)} s^k [f_\alpha(s^{-1/\alpha})]^{-1} \int_0^1 p^{n-1-k\alpha} f_\alpha((1-p)s^{-1/\alpha})dp \prod_{j=1}^k (1 -\alpha)_{n_j-1},
$$ 
hence by Bayes' rule 
$$
f_{S_\alpha, \gamma}(s|n_1, \dots, n_k)=\frac{ p_{\alpha, \gamma}(n_1, \dots, n_k| s^{-1/\alpha}) \gamma(s^{-1/\alpha})}{\int_0^\infty  p_{\alpha, \gamma}(n_1, \dots, n_k| s^{-1/\alpha}) \gamma(s^{-1/\alpha})ds }
$$
which simplifies to
$$
f_{S_{\alpha, h}}(s|K_n=k)=\frac{h(s^{-1/\alpha}) s^{k-1/\alpha -1} \int_0^1 p^{n-1-k\alpha} f_\alpha((1-p)s^{-1/\alpha})dp }{\int_0^\infty h(s^{-1/\alpha}) s^{k-1/\alpha -1} [\int_0^1 p^{n-1-k\alpha} f_\alpha((1-p)s^{-1/\alpha})dp] ds }
$$
and the result is proved. 
Notice that by definition of mixed $PK(\rho_\alpha, h \times f_\alpha)$  model, the general weights $V_{n,k, h}$ in the EPPF (\ref{eppf}) arise as follows
$$
V_{n,k, h}= \frac{\alpha^{k -1}}{\Gamma(n -k\alpha)} \int_0^\infty  h(s^{-1/\alpha}) s^{k -1/\alpha -1}  \int_0^1 p^{n-1-k\alpha} f_\alpha((1-p)s^{-1/\alpha}) dp ds
$$
hence the normalizing constant in formula (\ref{conddiv}) may be obtained through the following relationship (see also Ho et al. 2008, eq. (12))
\begin{equation}
\label{weiform}
\mathbb{E}_{n,k}^{\alpha} [h(S^{-1/\alpha})]=V_{n,k, h} \frac{\alpha^{1-k} \Gamma(n)}{\Gamma(k)}.
\end{equation}
{\bf Remark 2.} The result in Theorem 1. agrees with an analogous result for the conditional distribution $T| K_n=k$ (where $T$ is the random total sum of the ranked atoms of the r..p.m. $P$) first derived in an unpublished manuscript by Ho et al. (2008, cfr. Eq. (13)), that we received by one of those authors as a personal communication. Their result relies on the r.v. 
$$
R_{\alpha, (n,k)}= \frac{S_{\alpha, k\alpha}}{\beta({k\alpha, n-k\alpha})} 
$$
for $S_{\alpha, k\alpha}$ the polynomially tilted stable random variable with density
$$
f_{S_{\alpha, k\alpha}} (t)= \frac{\Gamma(k\alpha +1)}{\Gamma(k +1)} t^{ -k\alpha} f_\alpha(t). 
$$
It is an easy task to show that $[R_{\alpha, (n,k)}]^{-\alpha}= [Y_{\alpha, k} \times \beta(k\alpha, n- k\alpha)^{\alpha}]$. 
\\\\
In what follows we show how some distributional results for the conditional $\alpha$ diversity of specific Poisson-Kingman models already obtained in the literature may be easily derived by the general formula in Theorem 1.\\\\ 
{\bf Example 3.} [{\it Two-parameter Poisson-Dirichlet $(\alpha, \theta)$ partition models}] A first result for the conditional $\alpha$ diversity for the two-parameter Poisson-Dirichlet model (Pitman \& Yor, 1997) has been derived in Favaro et al. (2009) in view of Bayesian nonparametric posterior interval estimation for the number of new species in an additional sample in species sampling problems. Those authors rely on mimicking the original proof for the unconditional $\alpha$ diversity in Pitman (2006, Th. 3.8). To apply the general result it is enough to notice that the two-parameter Poisson-Dirichlet $(\alpha, \theta)$ for $\theta > -\alpha$ partition model is well-known to correspond (cfr. Pitman, 2003) to a mixed Poisson-Kingman model driven by the stable subordinator with mixing density 
$$
\gamma_{\alpha, \theta}(t)=  h(t) \times f_\alpha(t)= \frac{\Gamma(\theta+1)}{\Gamma(\theta/\alpha +1)} t^{-\theta} f_\alpha(t)
$$
with the weights in the Gibbs representation of the EPPF corresponding to
$$
V_{n,k}^{\alpha, \theta}=\alpha^{k-1}\frac{\Gamma(\theta/\alpha +k) \Gamma(\theta +1)}{\Gamma(\theta +n) \Gamma(\theta/\alpha +1)}. 
$$
Hence, by (\ref{weiform}), the denominator in (\ref{conddiv}) corresponds to
$$
\mathbb{E}_{n,k}^{\alpha}(h(Z^{-1/\alpha}))=\mathbb{E}_{n,k}^{\alpha}(Z^{\theta/\alpha})= \frac{\Gamma(\theta/\alpha +k) \Gamma(\theta +1)}{\Gamma(\theta +n) \Gamma(\theta/\alpha +1)} \frac{\Gamma(n)}{\Gamma(k)}
$$
and by Theorem 1. the conditional $\alpha$ diversity $Z_{n,k}^{\alpha, \theta}$ has density
\begin{equation}
\label{pdcond}
f_{n,k}^{\alpha, \theta}(z) =\frac{\Gamma(\theta +n)}{\Gamma(n-k\alpha) \Gamma(\theta/\alpha +k)} z^{\theta/\alpha +k -1-1/\alpha} \int_0^1 f_\alpha[(z w^{-\alpha})^{-1/\alpha}] (1 -w)^{n -k\alpha -1} dw
\end{equation}
which corresponds to the  $\theta/\alpha$ polynomial tilting of $\tilde{g}_{n,k}^{\alpha}$.  It is an easy task to verify that this is the density of the product of independent r.v.s $Z_{n,k}^{\alpha, \theta}=Y_{\alpha, \theta/\alpha +k} \times [\beta(\theta +k\alpha, n-k\alpha]^{\alpha}$,  
as already established in Cerquetti (2011) moving from a decomposition approach exploiting the {\it deletion of classes property} of this model (see Gnedin et al., 2009). In the next Proposition we prove the result actually agrees with the original result in Favaro et al. (2009) expressed in terms of the alternative scale mixture representation $Y_{\alpha, (\theta +n)/\alpha} \times \beta(\theta/\alpha +k,  n/\alpha -k)$. Notice that the equivalence in distribution between $R_{\alpha, (n,k)}$ and $S_{\alpha, n} [\beta_{k, n/\alpha -k}]^{-1/\alpha}$ was already stated without proof in Ho et al. (2008, cfr. Prop. 2.1, Eq. (11)). \\\\
{\bf Proposition 4.} Let $H=Y_1\times X$  for  $Y_1$ and $X$ independent r.v.s, $Y_1 \sim g_{\alpha, (\theta+n)}$ and $X\sim \beta(\theta/\alpha +k,  n/\alpha -k)$, then  the r.v. $Z_{n,k}^{\alpha, \theta}$ with density (\ref{pdcond}) and $H$ have the same characteristic function
%let $Z=Y_2*W^\alpha$ for $Y_{(\theta+k\alpha)/\alpha}$ and $W \sim Beta(\theta +k\alpha, n-k\alpha)$ then  
$$
G_{\alpha, \theta}^{n,k}(t)= \sum_{r \geq 0} \frac{(it)^r}{r!} \left(\frac{\theta +k\alpha}{\alpha}\right)_r
 \frac{1}{(\theta +n)_{r\alpha}}.
$$
\\
{\it Proof}:   First notice that by Proposition 2. in Cerquetti (2011), for $m \rightarrow \infty$,
$$
\mathbb{E}_{\alpha, \theta }\left(\frac{K_m^r}{m^{r\alpha}}\Big\vert  K_n=k\right) \stackrel{}{\longrightarrow} \left(\frac{\theta +k\alpha}{\alpha}\right)_r \frac{\Gamma (\theta +n)}{\Gamma(\theta +n +r\alpha)}.
$$
By the change of variable $zw^{-\alpha}=s$ (\ref{pdcond}) may be written as
$$
%\label{charac}
f_{n,k}^{\alpha, \theta}(z)=\frac{\Gamma(\theta +n)\Gamma(\theta+k\alpha+1)}{\Gamma(\theta +k\alpha)\Gamma(n -k\alpha) \Gamma((\theta  +k\alpha)/\alpha +1)}\frac 1\alpha z^{\theta/\alpha+k-1}\times
$$
$$\times \int_z^{\infty}\alpha^{-1}s^{-1/\alpha-1} f_\alpha(s^{-1/\alpha})\left( 1-(z/s)^{1/\alpha}\right)^{n-k\alpha-1}ds.
$$
Its characteristic function is given by 
$$
G_{n,k}^{\alpha, \theta}(t)=\frac{\Gamma(\theta +n)\Gamma(\theta+k\alpha+1)}{\Gamma(\theta +k\alpha)\Gamma(n -k\alpha) \Gamma((\theta  +k\alpha)/\alpha+1)} \frac 1\alpha \times $$
$$ \times \int_0^\infty \exp\{itz\} z^{\theta/\alpha+k-1}\int_z^{\infty}g_{\alpha}(s)\left( 1-(z/s)^{1/\alpha}\right)^{n-k\alpha-1}ds dz
$$
and may be rewritten as
$$
G_{n,k}^{\alpha, \theta}(t)=\frac{\Gamma(\theta+k\alpha+1)}{\Gamma((\theta  +k\alpha)/ \alpha+1)} 
\frac 1\alpha \int_z^{\infty}g_{\alpha}(s) \times
$$
$$ \times
\int_0^\infty \exp\{itz\}\frac{ \Gamma(\theta +n)}{ \Gamma(\theta +k\alpha)\Gamma(n -k\alpha)}z^{\theta/\alpha+k-1}
\left( 1-(z/s)^{1/\alpha}\right)^{n-k\alpha-1}dz ds \nonumber.
$$
By the change of variable $(z/s)^{1/\alpha}=y$, $z=y^{\alpha}s$, $dz=s\alpha y^{\alpha-1}dy$ yields
$$
=\frac{\Gamma(\theta+k\alpha+1)}{ \Gamma((\theta  +k\alpha)/\alpha+1)} 
\frac 1\alpha \int_0^{\infty} g_{\alpha}(s)
\times
$$
$$ \times \int_0^s e^{ity^{\alpha}s}\frac{\Gamma(\theta +n)}{ \Gamma(\theta +k\alpha)\Gamma(n -k\alpha)}(y^{\alpha}s)^{\theta/\alpha+k-1}
\left( 1- y\right)^{n-k\alpha-1}s \alpha y^{\alpha -1}dy ds
$$
and then reduces to
%\begin{equation}
%=\frac{\Gamma(\theta+k\alpha+1)}{ \Gamma(\frac{\theta  +k\alpha}{\alpha}+1)} 
%\int_0^{\infty} s g_{\alpha}(s)
%\int_0^s e^{ity^{\alpha}s}\frac{\Gamma(\theta +n)}{ \Gamma(\theta +k\alpha)\Gamma(n -k\alpha)}(y^{\alpha}s)^{\theta/\alpha+k-1}
%\left( 1- y\right)^{n-k\alpha-1}\alpha y^{\alpha -1}dy ds
%\end{equation}
$$
=\frac{\Gamma(\theta+k\alpha+1)}{ \Gamma((\theta  +k\alpha)/\alpha +1)} 
\int_0^{\infty} s^{\theta/\alpha +k} g_{\alpha}(s)
\int_0^1 e^{ity^{\alpha}s}\frac{\Gamma(\theta +n)}{ \Gamma(\theta +k\alpha)\Gamma(n -k\alpha)}(y)^{\theta+k\alpha-1}
\left( 1- y\right)^{n-k\alpha-1}dy ds.
$$
Exploiting the known characteristic function of $Y^\alpha$ for $Y \sim Beta(\theta +k\alpha, n -k\alpha)$ we can write
$$
\label{last}
=\frac{\Gamma(\theta+k\alpha+1)}{ \Gamma((\theta  +k\alpha)/\alpha +1)} \sum_{r=0}^{\infty} \frac{(it)^{r}}{r!}\frac{(\theta +k\alpha)_{r\alpha} }{(\theta +n)_{r\alpha}}  \int_0^{\infty} s^{\theta/\alpha +k +r} g_{\alpha}(s) ds
$$
and by (\ref{polmit})
\begin{equation}
\label{last}
=\sum_{r=0}^{\infty} \frac{(it)^{r}}{r!}\frac{(\theta +k\alpha)_{r\alpha }}{(\theta +n)_{r\alpha}} \frac{\Gamma(\theta+k\alpha+1)}{ \Gamma((\theta  +k\alpha)/\alpha +1)}  \frac{\Gamma((\theta +k\alpha +r\alpha)/\alpha +1)}{\Gamma(\theta +k\alpha +r\alpha +1)}.
\end{equation}
By the usual properties of Gamma function the last expression corresponds to
$$
=\sum_{r=0}^{\infty} \frac{(it)^{r}}{r!}\frac{\Gamma(\theta + k\alpha +r\alpha)\Gamma(\theta +n)}{\Gamma(\theta +k\alpha) \Gamma(\theta +n +r\alpha)}\frac{(\theta +k\alpha)\Gamma(\theta +k\alpha)}{\frac{\theta +k\alpha}{\alpha}\Gamma(\frac{\theta +k\alpha}{\alpha}) }\frac{\Gamma(\frac{\theta +k\alpha}{\alpha}+r) \frac{\theta +k\alpha +r\alpha}{\alpha}}{(\theta +k\alpha +r\alpha) \Gamma(\theta +k\alpha +r\alpha)}
$$
which simplifies to
\begin{equation}
\label{chfunc}
=\sum_{r=0}^{\infty} \frac{(it)^{r}}{r!}\left(\frac{\theta +k\alpha}{\alpha}\right)_r\frac{1}{(\theta +n)_{r\alpha}}
\end{equation}
and the conclusion follows by the result in Proposition 2 in Favaro et al. (2009) proving that (\ref{chfunc}) is the characteristic function of $H=Y_1 \times X$.\\\\
{\bf Example 5.} [{\it Generalized Gamma partition models}] As from Pitman (2003), generalized Gamma partitions models belong to the Poisson-Kingman family driven by the stable subordinator $PK(\rho_\alpha, \gamma)$ for a mixing density 
$$
\gamma_{\alpha, \lambda}(t)= \exp\{\psi_\alpha(\lambda) -\lambda t\} f_\alpha(t),
$$
where $\psi_\alpha(\lambda)=(2\lambda)^{\alpha}$ is the  Laplace exponent of $f_\alpha(\cdot)$. By an application of (\ref{prior}), after the reparametrization  $\lambda= \beta^{1/\alpha}/2$, the {\it unconditional} $\alpha$ diversity for this model (see also Cerquetti (2007) and Lijoi et al. (2007b)) is given by 
$$
\gamma_{\alpha, \beta}(s^{-1/\alpha})= \exp \left\{\beta -\frac 12 \left(\frac {\beta}{\alpha}\right)^{1/\alpha}\right\}f_\alpha(s^{-1/\alpha}) \alpha^{-1} s^{-1/\alpha -1}.
$$ 
To obtain the density of the conditional {\it posterior} $\alpha$ diversity it is enough to apply formula (\ref{conddiv}) which provides 
$$
f_{n,k}^{\exp, \alpha}(s)= \frac{\exp\left\{\beta -\frac 12\left(\frac{\beta}{s}\right)^{1/\alpha}\right\} \tilde {g}_{n,k}^{\alpha}(s)}{\mathbb{E}_{n,k}^{\alpha}\left[\exp\left\{\beta - \frac 12\left(\frac{\beta}{S}\right)^{1/\alpha}\right\}\right]}.
$$
The denominator arises by (\ref{weiform}) and the known expression for the  $V_{n,k}$ of $PK(\rho_\alpha, \gamma_{\alpha, \beta})$ models  as obtained in Eq. (6) in Cerquetti (2007),
\begin{equation}
\label{weingg}
V_{n,k}^{\alpha, \beta}= \frac{e^{\beta} 2^n \alpha^{k}}{\Gamma(n)} \int_{0}^{\infty} \lambda^{n-1} \frac{e^{-(\beta^{1/\alpha} +2\lambda)^{\alpha}}}{(\beta^{1/\alpha} +2\lambda)^{n-k\alpha}} d\lambda,
\end{equation}
which rewritten  in terms of incomplete Gamma functions by the change of variable $(\beta^{1/\alpha} + 2\lambda)^{\alpha}= x$, $d\lambda= (2\alpha)^{-1} x^{1/\alpha -1}dx$ yields 
$$
V_{n,k}^{\alpha, \beta}=\frac{e^{\beta} \alpha^{k-1}}{\Gamma(n)} \sum_{i=0}^{n-1} {n -1 \choose i} (-1)^{i}(\beta)^{i/\alpha} \Gamma(k - \frac{i}{\alpha}; \beta).
$$
By equation (\ref{weiform})
$$
{\mathbb{E}_{n,k}^{\alpha}\left[\exp\left\{\beta - \frac 12\left(\frac{\beta}{S}\right)^{1/\alpha}\right\}\right]}= \frac{e^{\beta} }{\Gamma(k)} \sum_{i=0}^{n-1} {n -1 \choose i} (-1)^{i}(\beta)^{i/\alpha} \Gamma(k - \frac{i}{\alpha}; \beta)
$$
hence by  Theorem 1. the conditional $\alpha$ diversity $S_{n,k}^{\alpha, \beta}$ for the generalized Gamma model has density 
$$
f_{n,k}^{\alpha, \beta}(s)=\frac{\Gamma(k) \exp (2^{-1}(\beta/s)^{1/\alpha}) \tilde{g}_{n,k}^{\alpha}(s)}{ \sum_{i=0}^{n-1} {n -1 \choose i} (-1)^{i}(\beta)^{i/\alpha} \Gamma(k - \frac{i}{\alpha}; \beta)}.
$$
The result agrees with Favaro et al. (2011, Th. 1.) due to the equivalence in distribution between 
$Y_{\alpha, n/\alpha} \times \beta (k, n/\alpha -k)$ and $Y_{\alpha, k} \times \beta(k\alpha, n- k\alpha)$ which arises specializing for $\theta=0$ the result in Proposition 4.

\section*{References}
\newcommand{\bibu}{\item \hskip-1.0cm}
\begin{list}{\ }{\setlength\leftmargin{1.0cm}}

\bibu \textsc{Cerquetti, A.} (2007) A note on Bayesian nonparametric priors derived from exponentially tilted Poisson-Kingman models. {\it Stat \& Prob Letters}, 77, 18, 1705--1711.

\bibu \textsc {Cerquetti, A.} (2009) A Generalized sequential construction of exchangeable Gibbs partitions with application. {\it Proceedings of S.Co. 2009, September 14-16, Milano, Italy. }

\bibu \textsc {Cerquetti, A.} (2011) A decomposition approach to Bayesian nonparametric estimation under two-parameter Poisson-Dirichlet priors. {\it Proceedings of ASMDA 2011 - Applied Stochastic Models and Data Analysis. Rome, June, 7-10 2011}

%\bibu  \textsc {Charalambides, C. A.} (2005) {\it Combinatorial Methods in Discrete Distributions}. Wiley, Hoboken NJ.

\bibu \textsc{Favaro, S., Lijoi, A., Mena, R.H. and Pr\"unster, I.} (2009) Bayesian non-parametric inference for species variety with a two-parameter Poisson-Dirichlet process prior. {\it JRSS-B}, 71, 993-1008.

\bibu \textsc{Favaro, S., Lijoi, A. and Pr\"unster, I.} (2011) Asymptotics for a Bayesian nonparametric estimator of species variety. {\it Bernoulli} (to appear).

\bibu \textsc{Gnedin, A., Haulk, S. and Pitman, J.} (2009) {Characterizations of exchangeable partitions and random discrete distributions by deletion properties}. {\sf http://arxiv.org/abs/0909.3642}

\bibu \textsc{Gnedin, A. and Pitman, J. } (2006) {Exchangeable Gibbs partitions  and Stirling triangles.} {\it Journal of Mathematical Sciences}, 138, 3, 5674--5685. 

%\bibu \textsc{Hansen, B. and Pitman, J.} (2000) Prediction rules for exchangeable sequences related to species sampling. {\it Statistics \& Probability Letters}, 46, 251--256.

\bibu \textsc{Ho, M-W, James, L.F. and Lau, J.W.} (2008) Explicit Gibbs Chinese Restaurant Process priors. {\it Unpublished manuscript}

%\bibu \textsc{Ho, M-W, James, L.F. and Lau, J.W.} (2007) Gibbs partitions (EPPF's) derived from a stable subordinator are Fox H - And Meijer G - Transforms. arXiv:0708.0619v2 [math.PR]

%\bibu \textsc{Hsu, L. C, \& Shiue, P. J.} (1998) A unified approach to generalized Stirling numbers. {\it Adv. Appl. Math.}, 20, 366-384.

%\bibu \textsc{Johnson, N.L. \& Kotz, S.} (1977) {\it Urn models and their application}. Wiley \& Sons.

%\bibu \textsc{Johnson, N. L. \& Kotz, S.} (2005) {\it Univariate Discrete Distributions}. 3rd Edition, Wiley \& Sons.

\bibu \textsc{Lijoi, A., Mena, R.H. and Pr\"unster, I.} (2007a) Controlling the reinforcement in Bayesian non-parametric mixture models. {\it JRSS-B}, 69, 4, 715--740. 

\bibu \textsc{Lijoi, A., Mena, R.H. and Pr\"unster, I.} (2007b) Bayesian nonparametric estimation of the probability of discovering new species.  {\it Biometrika}, 94, 769--786.

\bibu \textsc{Lijoi, A., Pr\"unster, I. and Walker, S.G.} (2008) Bayesian nonparametric estimator derived from conditional Gibbs structures. {\it Annals of Applied Probability}, 18, 1519--1547.

%\bibu \textsc {Normand, J.M.} (2004) Calculation of some determinants using the $s$-shifted factorial. {\it J. Phys. A: Math. Gen.} 37, 5737-5762.

%\bibu \textsc{Perman, M., Pitman, J, \& Yor, M.} (1992) Size-biased sampling of Poisson point processes and excursions. {\it Probab. Th. Rel. Fields}, 92, 21--39.

%\bibu \textsc{Pitman, J.} (1995) Exchangeable and partially exchangeable random partitions. {\it Probab. Th. Rel. Fields}, 102: 145-158.

%\bibu \textsc{Pitman, J.} (1996a) Some developments of the Blackwell-MacQueen urn scheme. In T.S. Ferguson, Shapley L.S., and MacQueen J.B., editors, {\it Statistics, Probability and Game Theory}, volume 30 of {\it IMS Lecture Notes-Monograph Series}, pages 245--267. Institute of Mathematical Statistics, Hayward, CA.

%\bibu \textsc{Pitman, J.} (1996b) Notes on the two parameter generalization of Ewens random partition structure. {\it Manuscript} University of California, Berkeley. Unpublished.

\bibu \textsc{Pitman, J.} (2003) {Poisson-Kingman partitions}. In D.R. Goldstein, editor, {\it Science and Statistics: A Festschrift for Terry Speed}, volume 40 of Lecture Notes-Monograph Series, pages 1--34. Institute of Mathematical Statistics, Hayward, California.

\bibu \textsc{Pitman, J.} (2006) {\it Combinatorial Stochastic Processes}. Ecole d'Et\'e de Probabilit\'e de Saint-Flour XXXII - 2002. Lecture Notes in Mathematics N. 1875, Springer.

\bibu \textsc{Pitman, J. and Yor, M.} (1997) The two-parameter Poisson-Dirichlet distribution derived from a stable subordinator. {\it Ann. Probab.}, 25, 855--900.

%\bibu \textsc{Toscano, L.} (1939) Numeri di Stirling generalizzati operatori differenziali e polinomi ipergeometrici. {\it Comm. Pontificia Academica Scient. } 3:721-757.

%\bibu \textsc{Yamato, H. and Sibuya, M.} (2000) Moments of some statistics of Pitman sampling formula. {\it Bull. Inform. Cybernet.}, 32 1--10.

\end{list}
\end{document}